# MANY-AGENT INTERACTION IN THE MODEL OF LABOUR FORCE TRAINING

[1, a)]**IRINA ZAITSEVA**, [2, b)]**OLEG MALAFEYEV**, [2, c)]**SERGEI STREKOPYTOV**, [1, d)]**ANNA ERMAKOVA**, [1, e)] **DMITRY SHLAEV**

[1]Stavropol State Agrarian University, 12, Zootechnicheskiy alley, Stavropol, Russia, 355017

[2]Saint-Petersburg State University, 7-9, University Emb., Saint-Petersburg, Russia, 199034

E-mail: [a)]zirinazirina2015@yandex.ru, [b)]malafeyevoa@mail.ru, [c)]SergeiStrekopytov@mail.ru, [d)]shl-dmitrij@yandex.ru, [e)]dannar@list.ru

**ABSTRACT**

The continuous and discrete models of labour force training are being built. The application of the results from the theory of differential games and dynamic programming allows presenting the optimal strategies of labour force training that can be calculated.

**Keywords:** *Continuous game-theoretical model, discrete game-theoretical model, labour force training.*

## 1. INTRODUCTION

The economic demand on "educated" (high-skill) employees has been grown recently. The expected discounted salaries of high-skill employees have grown. Thus, the investments of time, effort and money into education have been attracting more and more citizens. The "educated employees" are supposed to be competent workers that have all the required skills and knowledge to accomplish their responsibilities.

The education can reduce unemployment all over the world and Russia is not the exception. Education is an investment project with additional positive or negative benefit. The educational service is characterized by immateriality and variability; it's difficult for consumers to estimate its quality. It's also hard to compare this service with analogous services. These facts stipulated the complexity of labour force training research.

The problem of labour force training research and necessity to enlarge the employees' level of education can be solved via economic-mathematical methods.

## 2. MODELING OF LABOR MARKETS AND EDUCATIONAL SERVICES

Relying on one of the most important elements of the state policy in the field of optimizing the costs of education, which is connected with the formation and development of the education system - the labor market - the economy. This structural model reflects the picture of the dynamics of labor resources: at the input - the birth rate vector, at the output - the gross regional product (GRP). From the position of macroeconomics in this scheme, a person acts as an impersonal labor resource necessary for the production of goods and services [1].

The need for specialists with higher education for the economy can be calculated in a standard way, based on a change in certain factors shown in Figure 1 [2].

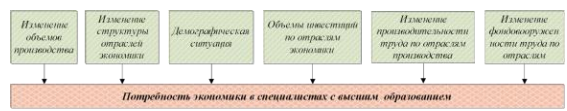

*Figure 1* - Factors determining the need for specialists with higher education

Under the influence of these factors, the region's need for cadres with a higher education in the year **t** on specialty **i** is formed.

$$R_t^i = Re_t^i + Ro_t^i \qquad (1)$$

Where $R_t^i$ – the need of the state in the year t in staff with a higher education in the specialty ***i***;

$Re_t^i$ – need of the economy in a year t in specialists with higher education in the specialty ***i***;

$Ro_t^i$ – need of society in a year t in obtaining higher and postgraduate education in the specialty ***i***.

Participants in the educational process, between which the material relationship takes place, are:

students, organizations providing educational services; persons and organizations that pay for these services (in the role of which are the societies and the labor market).

Before solving the problem of synthesis of SDS (complex dynamic systems), it is necessary to solve the problem of their analysis, which consists in decomposition of a complex system into subsystems and establishment of connections between them.

The representation of the labor market (RT) as a manager makes it possible, taking into account the principles of system analysis, to detail the interconnected modules with the subsequent possibility of building a low-level model with their integration into the global complex model.

Individuals demand for higher education as an investment asset, which also has some consumer value. Having studied, the graduate offers his labor in the labor market of educated specialists at the equilibrium rate for his specialty and works until the end of his career. The general equilibrium in the labor market of specialists and higher education determines the wages of skilled workers, their share in the employed population, the share of students in the population from 17 to 21 years, and the equilibrium number of HEIs and their optimal size. The salary of unskilled workers is exogenous and constant due to the presence of an unlimited number of such workers.

It is also assumed that the rate of change in output in the economy is significantly lower than the speed described in the model of processes that we have constructed, so we assume that aggregated output is exogenous and fixed. As a result, changes in technology are caused only by changes in the demand for skilled and unskilled labor, and changes in output occur after equilibrium is established in the education and labor markets.

The economy consists of two sectors: the national economy (including the service sector, except for higher education) and the system of higher education.

Retaining the notation used in the work of Shevchuk D.V [1] and supplementing them with indices, we have. All population L divided into educated workers C , teachers E , students S and uneducated workers N , here the equality N + C + E + S = L. Thus, there are three decision-makers in the economy: individuals, the labor market and universities. It is advisable to consider each of them separately.

Individuals live two unequal periods: at the beginning of the first one everyone decides to work for him as an uneducated worker or to study and does it all the time; in the second he works in accordance with his level of education.

When deciding on education, the individual tries to maximize his life-long discounted income, which is the sum of the terms of the geometric progression with the denominator β and the first member $w_N$ , equal to some equilibrium salary of an uneducated person. From the theory of analysis of financial rents, it follows that the income is determined by the formula:

$$\sum_{s=0}^{T} \beta^s w_N = w_N \frac{1-\beta^T}{1-\beta}.$$

This stream is compared with the flow of an educated person who will pay for the year of study for $i$ specialty price p($i$) and immediately after graduating from the university will begin to receive the equilibrium annual salary of an educated specialist $w_h(i)$ for simplicity in the form of an annuity:

$$-\sum_{s=0}^{x-1} \beta^s p(i) + \beta^x \sum_{s=0}^{T-x} \beta^s w_h(i) = -p(i)\frac{1-\beta^{x-1}}{1-\beta} + \beta^x w_h(i)\frac{1-\beta^{T-x}}{1-\beta}.$$

If $x$, the number of years of study is small, for example, five, then for the same reasons as in the case of an uneducated individual, one can omit the term $\beta^{T-x}$, then $\frac{1-\beta^{T-x}}{1-\beta} = 1$. In real life, of course, the salary of a specialist increases with the experience and the main income falls on the second half of the career. This means that the change in the denominator of the progression β will affect the discounted income more than in the case of an annuity, but the quality (sign) of changes will remain the same.

Thus, when making a positive decision about training, the difference between salaries should be large enough to cover the cost of training, taking into account the degree of impatience of the individual.

The next element of the market of educational services, which actually produces these services, is a higher educational institution. The production of educated labor is carried out by technology with increasing returns to scale: for the production of a certain number of specialists $h$ need to spend a year $F + ch$ units of skilled work of a specialty, which the university itself produces. Constant costs are generally considered to the costs of maintaining the general management apparatus of the university; its positions are unshakable, especially in state institutions. The vacant premises can be rented, and ordinary teachers can change the load, i.e. costs can be attributed to variable costs $c$.

Every year the university produces $h_{\Sigma}(i)$ specialists $i$, part of which replaces the retired labor in the national economy, and the other - in education. If we assume that generations have the same number, then every year in the national economy it is eliminated $\frac{C}{T-x}$ skilled labor, and in education – $\frac{E}{T-x}$. Consequently, in sum, all universities should produce each year $\frac{E+C}{T-x}$ graduates.

The university's desire to maximize extrabudgetary income can be characterized by a function:

$$\pi(i) = x\big[p(i)h_{\Sigma}(i) - w_h(i)(F + ch_{\Sigma}(i))\big].$$

(2)

Note that the assumption that the university pays specialists as much as the national economy makes sense for the balance. The University makes a decision when assigning a certain price for educational services and the annual number of students. The price, as is known from the individual's task, must satisfy the condition $\alpha\beta^x w_h(i) - w_N \geq p(i)(1 - \beta^{x-1})$. In equilibrium, one can expect that the university will set the next price for an educational service:

$$p(i) = \frac{\alpha\beta^x w_h(i) - w_N}{1 - \beta^{x-1}}$$

(3)

since at this price you can recruit any number of students, individuals still learn or not. When the price decreases, the number of applicants will also be any within the population of the corresponding age, and the profit will decrease [1].

Since the price of training are determined by the difference in wages between skilled and unskilled workers, the balance in the labor market will affect the university decision for graduates, where their salary is determined. In other words, the university learns the current salary levels on the market and sets the price according to (3), and the output will be formed according to the maximization of profit condition.

The condition for the existence of education can be rewritten, taking into account (3) as a condition for the difference in wages of educated and uneducated people:

$$w_h(i)(\alpha\beta^x - c(1 - \beta^{x-1})) > w_N.$$

The cost of the discounted flow from the salary of a qualified person through x years minus the cost of the flow of education costs should be higher than the cost of flow from the salary of an unskilled person who is starting to work now. Note that the inequality in income between skilled and unskilled people is the main condition for the existence of higher education.

Then the specialist gets to the RT, composed of three main components, which was shown in Figure 2:

- relationships between employees and employers;
- relations between subjects of the Republic of Tajikistan and representatives
- (trade unions, employers' associations, employment services);
- relations between the subjects of the Republic of Tatarstan and the state.

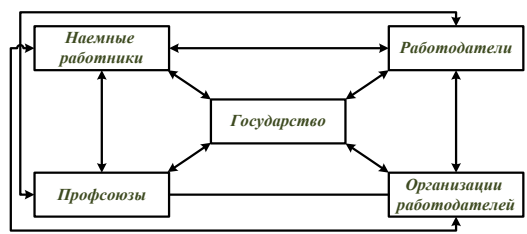

*Figure 2* - The system of relations in the labor market

The main components of the labor market are: the aggregate supply (P), covering all hired labor; and aggregate demand (C) as the general need of the economy for wage labor. They constitute the aggregate labor market.

Thus, the part that is formed by crossing the aggregate demand and the aggregate supply is called the satisfied demand for labor (SDL). Non-intersecting parts correspond to the current market:

CM = AM – SDL,

where AM is the aggregate labor market; CM - current market.

The current labor market is formed due to the natural and mechanical movement of labor and workplaces (the introduction of new and retirement of old ones).

Thus, it can be concluded that the more knowledgeable the entrants about the status of RT and the prestige of the university in this market, the share of the satisfied labor demand (SDL) in the total labor market is greater. Consequently, the unemployment rate will fall.

Despite the strict, at first glance, mathematical calculations, the origins in the basics of financial mathematics, the methodology only a qualitative analysis of a limited number of parties and the field of employment in connection with the inability to specify precise quantitative values of model parameters. The foregoing version modeling is mainly qualitative, explanatory character.

## THE CONTINUOUS MODEL OF LABOUR FORCE TRAINING

Let's observe the continuous model. We shall design $\tau$ as the time - $\in R_1^+ = [0, \infty)$, as $t \in R_1^+$ - the rate of labour force training process. $x$ - shall be designed as the level of labour force education, $x \in R_1^+$. The rate of education t is considered to be the control parameter that is managed by the employer. Let's denote as $\alpha$ the vector of k-dimensional Euclidean space $R^k$, identifying the external conditions that influence on the labour force training. This parameter is controlled. It is supposed not to be predicted. The dynamics of the labour force training process is defined by the differential equation:

$$\dot{x} = f(x, t, \alpha),$$

where $x(0) = x_0$ is the initial level of labour force education, $f$ is the continuous function satisfying a Lipchitz condition with respect to $x$:

$$|f(x, t, \alpha) - f(x', t, \alpha)| \leq K|x - x'|$$

at any $x, x' \in R_1^+, t \in T, \alpha \in Q$, where k-positive constant, $T$, $Q$ are variation ranges of control parameters $t$, $\alpha$ correspondently;

- a condition of solution time length:

$$|f(x, t, \alpha)| \leq M + N(x)$$

at any $x \in R_1^+$, $t \in T$, $\alpha \in Q$, where $M$, $N$ – positive constants;

- a condition of vectograms convexity:

Sets $f(x, T, \alpha) = \{y = f(x, t, \alpha)| t \in T\}$, sets $f(x, t, Q) = \{y = f(x, t, \alpha)| \alpha \in Q\}$ are convexed.

The measurable function $t(\tau)$, $t: R_1^+ \to T$ is called the admissible operator control. The measurable function $\alpha(\tau)$, $\alpha: [0, I^*] \to Q_-$ is called the admissible control of "uncontrolled variables".

The process path corresponding to the control couple $(t, \alpha)$ is called the absolute continuous function $x: [0, I^*] \to R_1^+$, is such as that almost at every $\tau \in [0, I^*]$ the relation is fulfilled.

$$\dot{x}(\tau) = f(x(\tau), t(\tau), \alpha(\tau)).$$

The process is supposed at interval $[0, I^*]$. Such a process is being chosen that at the time till the moment $I^*$ the path comes up to the set $\tilde{X} \subset R_1^+$, defined for the task in advance (game terminal set). The function $e(t, \alpha, x, \tau)$ characterizing training costs is also defined in the task. If $t(\tau)$, $\alpha(\tau)$ are controlled operator functions and "uncontrolled variables" at interval $[0, I^*]$, and $x(\tau)$ – is the corresponding process path, $\int_0^{I^*} e(t(\tau), \alpha(\tau), x(\tau), \tau) \alpha\tau = E(t, \alpha, x(t, \alpha))$ are the training costs. Let's define two tasks in the model framework:

1) using the control $t(\tau)$ (the rate of education) to bring the labour force level of education from the initial state to the final state $x(\tilde{\tau}) \in \tilde{X}$ with minimal costs at any external conditions-controls of "uncontrolled variables";

2) likewise, choosing the control $t(\tau)$ to bring the level of education to the defined condition $\tilde{X}$ at minimal time and at any "uncontrolled variables".

To be definite let's consider the first task. According to the meaning of the task it is natural to consider the upper (majorant) game $\bar{\Gamma}$ [3,5]. The piecewise programming strategies are defined in this game by the following. The strategy of the minimizing operator $\bar{T}$ is the couple $\bar{T} = (\sigma', \varphi_{\sigma'})$, where $\sigma'$ is the final partition of the interval $[0, I^*]$, $\varphi_{\sigma'}$ is the lower strategy corresponding to the partition $\sigma'$. This function assigns the employer's information at the moment $\tau_i \in \sigma'$ to the measurable control $t_i(\tau)$ at the interval $[\tau_i, \tau_{i+1}]$.

The strategy of uncontrolled variables $\tilde{\alpha}$ is the population $\{\psi_\sigma\}_{\sigma \in \Sigma}$, where $\Sigma$ is the set of all the final partitions of the interval $[0, I^*]$, $\psi_\sigma$ is the upper strategy of uncontrolled variables in the multistage upper game $\bar{\Gamma}^\sigma$ corresponding to the partition $\sigma'$. This function assigns the information of uncontrolled variables at the moments $t_i \in \sigma'$ to the measured control $\alpha_i(\tau)$ at the interval $[\tau_i, \tau_{i+1}]$. Let's remember that in the upper game $\bar{\Gamma}^\sigma$ the employer at any period of time $t_k \in \sigma$ knows about the state of the process at the particular period for uncontrolled variables at the moments $t_k \in \sigma$. Besides, the employer's control selected by him at the interval $[t_k, t_{k+1}]$ is also known. According to the couple strategies $(\bar{T}, \bar{\alpha})$ the unique set of control $t(\tau), \alpha(\tau)$ at the whole interval of the game $[0, I^*]$ and the corresponding process path $x(\tau) = \chi(\bar{T}, \alpha)$ are built. It is known that [3] there is ε-saddle point for any $\varepsilon > 0$ in the game $\bar{\Gamma}$. It means such a strategy couple $\bar{T}^\varepsilon, \bar{\alpha}^\varepsilon$, that at any $\bar{T} \in T, \bar{\alpha} \in Q$ the inequation is fulfilled.

$$\mathrm{E}(\bar{T}^\varepsilon, \bar{\alpha}, \chi(\bar{T}^\varepsilon, \bar{\alpha})) + \varepsilon \leq \mathrm{E}(\bar{T}^\varepsilon, \bar{\alpha}^\varepsilon, \chi(\bar{T}^\varepsilon, \bar{\alpha}^\varepsilon)) \leq \mathrm{E}(\bar{T}, \bar{\alpha}^\varepsilon, \chi(\bar{T}, \bar{\alpha})) - \varepsilon$$

The value $(\bar{T}^\varepsilon, \bar{\alpha}^\varepsilon, \chi(\bar{T}^\varepsilon, \bar{\alpha}^\varepsilon)) - \varepsilon$ is the final amount of material costs that could be spent on the labour force training at any uncontrolled variables. Thus, $E(\tau)$ is the bending in the situation $\bar{T}^\varepsilon \bar{\alpha}^\varepsilon$.

## 3. THE DISCRETE MODEL OF LABOUR FORCE TRAINING

It is known that the labour force must take the training in a fixed sequence in $n$ directions. At i-step $(i = \overline{1, n})$ depending on the rate of labour force education chosen by the employer that should be constant at this step and the vector of external conditions, $e_i(t_i, \alpha_i, x_{i-1})$ material costs are spent on the training. During the whole period of training at $n$ steps-operations the following amount of material costs is spent

$$\sum_{i=1}^{n} e_i(t_i, \alpha_i, 0x_{i-1}) = \varepsilon_n(t_1, \dots, t_n, \alpha_1, \dots, \alpha_n) = \varepsilon_n(t^n, \alpha^n)$$

The limits of possible variables in the rate of training are defined at every step $i$ and the vector of external conditions:

$$t_i \in T_i = [t_i^1, t_i^2]; \alpha_i \in Q_i = [\alpha_i^1, \alpha_i^2]$$

As in the 1st section, it is required to choose such a control sequence $t = \{t_i\}$, that can assure the least material costs at any possible external conditions $\alpha = \{\alpha_i\}$. It should be noticed that the labour force continue to be in i-state during the long period of time but the external conditions can change. So, in this case the set of states for the external conditions vector presents the set of functions at a time limit. It is natural to consider this set of functions as a compact space leaving the labeling Q.

It is supposed that the number of steps $n$ is such that the set of required level of education $\tilde{X}$ is achieved at any control $\{t_i\}$ and any functions of external conditions $\{\alpha_i\}$. We shall define by $F_n(x_0)$ the minimal assured material costs during the training process at $n$ steps, initial education level $x_0$ and any external conditions:

$$F_n(x_0) = \min_{t} \max_{\alpha} e_n(t_\alpha^n, \alpha^n, x^n)$$

It is easy to check that the method of dynamic programming is applied here. This method allows comprising recursion relation binding the values of the Bellman's function at two consequently connected values of n parameter

$$F_n(x_0) = \min_{t_1} \max_{\alpha_1}\{e_1(t_1, \alpha_1, x_0) - F_{n-1}(x_1)\},$$
$$x_1 = x_0 + \Delta t_1 f(x_0, t_1, \alpha_1).$$

## 4. CONCLUSION

The numerical evaluations done by PC for a discrete model at $n = 3$ verified theoretical calculations. The modeling of the labour force training becomes essential for the research of labour force training patterns, for the estimation of how the level of education influences on the economic development, for designing the methods and procedures of the labour force management. The content of the training process and use of labour force make urgent the research of the training peculiarities under the conditions of the market economy in Russia and in a particular region.

**ACKNOWLEDGEMENT**

The work is partly supported by work RFBR No. 18-01-00796.